%% file: saitoclean2.tex
\documentclass{article}
\usepackage{spconf,amsmath,amssymb,amsthm,euscript,fourier}
\usepackage{graphicx,version}
\usepackage[pdftex]{hyperref}
\excludeversion{hide}
\graphicspath{{../../../../src/julia/julia_test/}}
\include{amsmath_def}

\include{symbols}
\newcommand{\Path}{\operatorname{Path}}
\newcommand{\len}{\operatorname{length}}
\newcommand{\dtilde}[1]{\tilde{\tilde{#1}}}
\renewcommand{\ndim}{d}

\title{How can we naturally order and organize graph Laplacian eigenvectors?}
\name{Naoki Saito\thanks{This research was partially supported by the NSF grant
DMS-1418779 and the ONR grant N00014-16-1-2255; the author also thanks Prof.\ Qinglan Xia of UC Davis for the fruitful discussion and the \texttt{Plots.jl} community for answering his many graphics questions.}}
\address{University of California, Davis\\
  Department of Mathematics\\
  Davis, CA 95616 USA}

\begin{document}
\maketitle
\begin{abstract}
\input{abs2}
\end{abstract}
\section{Introduction}
\label{sec:intro}
For the theory and practice of discrete wavelets on regular lattices
in $\Rf^\ndim$, Fourier analysis has played a significant role.
Hence, when attempting to develop wavelet theory for graphs and networks,
some researchers have used graph Laplacian eigenvalues and eigenvectors in place
of the frequencies and complex exponentials, respectively; see, e.g.,
the Spectral Graph Wavelet Transform of Hammond et al.~\cite{HAMMOND-VANDERGHEYNST-GRIBONVAL}.
While tempting to do so, this viewpoint/strategy has several fundamental
problems.
One of them is the intricate relationship between the frequencies
and the Laplacian eigenvalues.
For undirected and unweighted paths (or cycles), the Laplacian eigenvectors are
the discrete cosine (or Fourier) basis vectors and the corresponding eigenvalues
are monotonically increasing functions of the frequency as discussed in \cite{SAITO-WOEI-DENDRITES, SAITO-WOEI-KOKYUROKU, WHY4-LAA, SAITO-BJSIAM-ENGLISH, IRION-SAITO-SPIE} among others.
Consequently on those simple graphs, one can precisely develop the classical
wavelets using the Littlewood-Paley theory \cite[Sec.~2.4]{Meyer-AppBook2}.
However, as soon as a graph becomes even slightly more complicated (e.g., a
discretized rectangle in 2D), the situation completely changes: we cannot
view the eigenvalues as a simple monotonic function of frequency anymore
 \cite{SAITO-BJSIAM-ENGLISH, IRION-SAITO-SPIE}.
Hence, a fundamental question is how to \emph{order} and \emph{organize}
Laplacian eigenvectors without using the eigenvalues, i.e., how to create
a \emph{dual} domain of a given graph.
In this paper, we investigate this important problem further and propose a new
method to order and organize those eigenvectors using the idea from the
Ramified Optimal Transport (ROT) theory \cite{XIA-REVIEW} to measure
``natural'' distances between eigenvectors followed by embedding them 
into a low-dimensional Euclidean domain.

\section{Problems with Ordering Eigenvectors According to the Corresponding Eigenvalues}
\label{sec:motivation}
Let $G=G(V, E)$ be a graph with its vertex set $V$ and edge set $E$.
It is a well-known fact by now \cite{SAITO-WOEI-DENDRITES, SAITO-WOEI-KOKYUROKU, WHY4-LAA, SAITO-BJSIAM-ENGLISH, IRION-SAITO-SPIE}
that if the graph $G$ is a path with $n$ nodes,
i.e., $G=P_n$, then the eigenvectors of the \emph{graph Laplacian matrix}
$L(P_n) \define D(P_n) - A(P_n)$, where $D(P_n)$, $A(P_n)$ are its degree and
adjacency matrices, respectively, are exactly the \emph{DCT Type II} basis
vectors (used for the JPEG standard) while those of the
\emph{symmetrically-normalized graph Laplacian matrix}
$L_\mathrm{sym}(P_n) \define D(P_n)^{-\half}\, L(P_n) \, D(P_n)^{-\half}$ are
the \emph{DCT Type I} basis vectors.
In fact, the eigenpairs of $L(P_n)$ are:
\begin{equation}
\label{eqn:eigen1}
\lambda_{k;n} \define 4 \sin^2\(\frac{\pi k}{2n}\) , \enskip
\phi_{k;n}[x] \define a_{k; n} \cos\(\frac{\pi k}{n}\(x+\half\)\),
\end{equation}
where $k, x = 0:n-1$, and $a_{k; n}$ is a normalization constant to have
$\| \bphi_{k;n} \|_2 = 1$. It is clear that the eigenvalue is a monotonically
increasing function of the \emph{frequency}, which is the eigenvalue index $k$
divided by $2$ in this case.

\begin{figure}[h]
\centering{\includegraphics[width=0.475\textwidth]{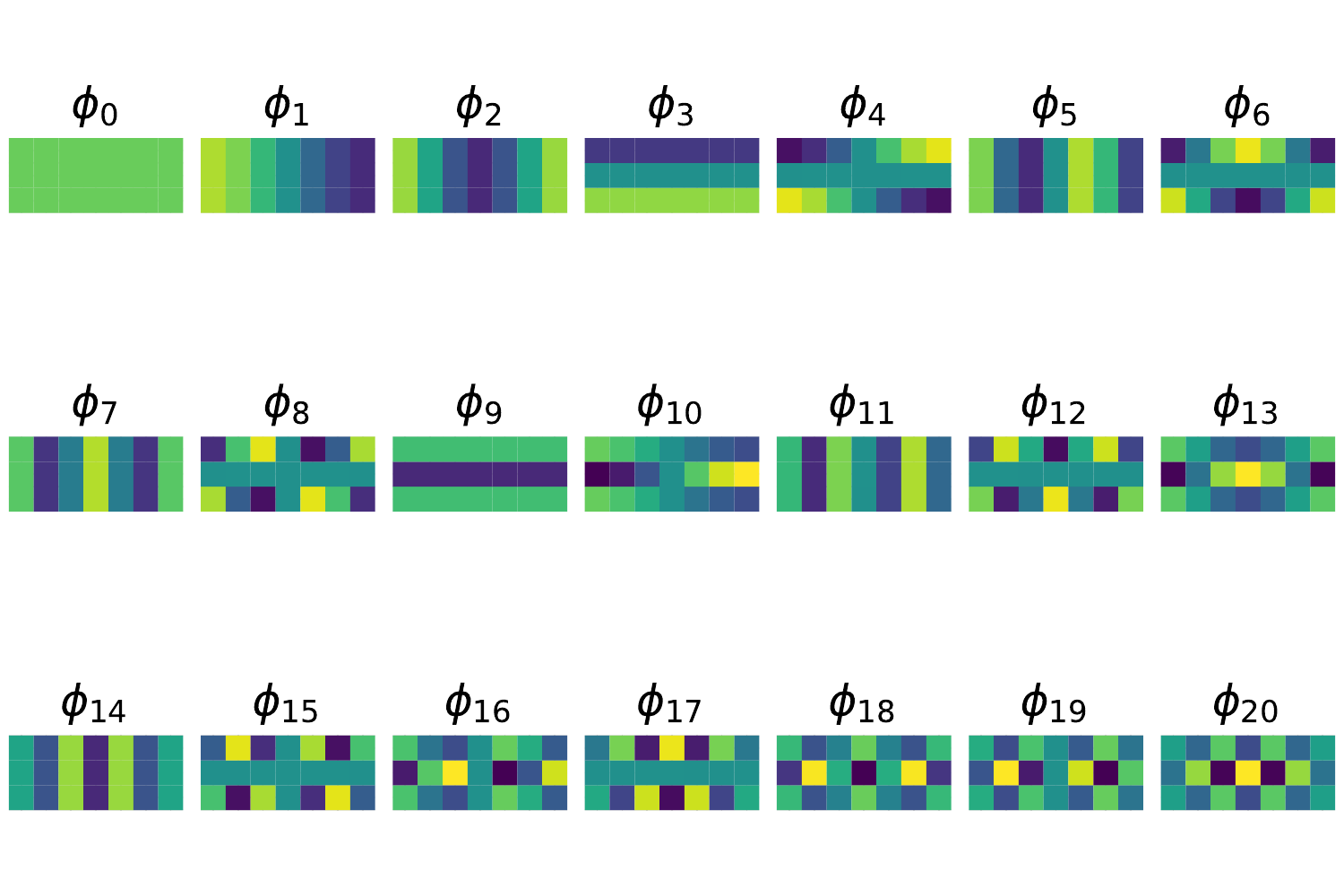}}
\vspace{-1.5em}
\caption{Laplacian eigenvectors of $P_7 \times P_3$ sequentially ordered in terms
of nondecreasing eigenvalues.}
\label{fig:grid7x3evorder}
\end{figure}

As soon as a graph becomes even slightly more complicated than
unweighted and undirected paths/cycles, however,
the situation completely changes: we cannot view the eigenvalues as a simple
monotonic function of frequency anymore.
For example, consider a thin rectangle in $\Rf^2$, and suppose that this
rectangle is discretized as $P_m \times P_n$ ($m > n > 1$). The Laplacian
eigenpairs of this graph can be easily derived from Eq.~\eqref{eqn:eigen1} as:
\begin{align*}
\lambda_k &= \lambda_{(k_x, k_y)} \define \lambda_{k_x; m} + \lambda_{k_y; n} \\
\phi_k[x,y] &= \varphi_{k_x,k_y}[x,y] \define \phi_{k_x; m}[x] \cdot \phi_{k_y; n}[y]
\end{align*}
where $k=0:mn-1$; $x, k_x=0:m-1$; and $y, k_y=0:n-1$.
As always, let $\{\lambda_k\}_{k=0:mn-1}$ be ordered in the
nondecreasing manner. Fig.~\ref{fig:grid7x3evorder} shows the corresponding
eigenvectors ordered in this manner (with $m=7$, $n=3$).
Note that the layout of $3 \times 7$ grid of subplots is 
for the page saving purpose: the layout of $1 \times 21$ grid of subplots would
be more natural if we use only the eigenvalue size for eigenvector ordering.
In this case, the smallest eigenvalue is still
$\lambda_0=\lambda_{(0,0)}=0$, and the corresponding eigenvector is constant.  
\vspace{-1em}
The second smallest eigenvalue $\lambda_1$ is
$\lambda_{(1,0)}=4\sin^2(\pi/2m)$, since $\pi/2m < \pi/2n$, and its
eigenvector has half oscillation (i.e., half period) in the $x$-direction.  
But, how about $\lambda_2$?
Even for such a simple situation
there are two possibilities for $\lambda_2$, depending on $m$ and $n$.
If $m > 2n$, then $\lambda_2=\lambda_{(2,0)} < \lambda_{(0,1)}$.
On the other hand, if $n < m < 2n$,
then $\lambda_2=\lambda_{(0,1)} < \lambda_{(2,0)}$.
More generally, if $Kn < m < (K+1)n$ for some $K \in \N$, then 
$\lambda_k=\lambda_{(k,0)}=4\sin^2(k\pi/2m)$ for $k=0, \dots, K$.
Yet we have $\lambda_{K+1}=\lambda_{(0,1)}=4\sin^2(\pi/2n)$ and
$\lambda_{K+2}$ is equal to either $\lambda_{(K+1,0)}=4\sin^2((K+1)\pi/2m)$
or $\lambda_{(1,1)}=4[\sin^2(\pi/2m) + \sin^2(\pi/2n)]$ 
depending on $m$ and $n$.
In Fig.~\ref{fig:grid7x3evorder}, one can see this behavior with $K=2$.
Clearly, the mapping between $k$ and $(k_x,k_y)$ is quite nontrivial.
Notice that $\phi_{(k,0)}$ has $k/2$ oscillations in the $x$-direction
whereas $\phi_{(0,1)}$ has only half oscillation in the $y$-direction.
In other words, all of a sudden the eigenvalue of a completely different
type of oscillation \emph{sneaks into} the eigenvalue sequence.
Hence, on a general graph, by simply looking at its
Laplacian eigenvalue sequence $\{\lambda_k\}_{k=0, 1, \dots}$,
it is \emph{almost impossible to organize the eigenvectors into physically
meaningful dyadic blocks and follow the Littlewood-Paley approach} unless
the underlying graph is of very simple nature, e.g., $P_n$ or $C_n$ (the cycle
consisting of $n$ nodes).
Therefore, for complicated graphs, wavelet construction methods that rely on
the Littlewood-Paley theory by viewing the Laplacian eigenvalues as the
frequencies, such as the spectral graph wavelet transform
\cite{HAMMOND-VANDERGHEYNST-GRIBONVAL}, may face unexpected problems.
In fact, not only the wavelet construction methods but also any procedures
and applications having that viewpoint would become problematic on general graphs.

\begin{figure}[h]
\centering{\includegraphics[width=0.475\textwidth]{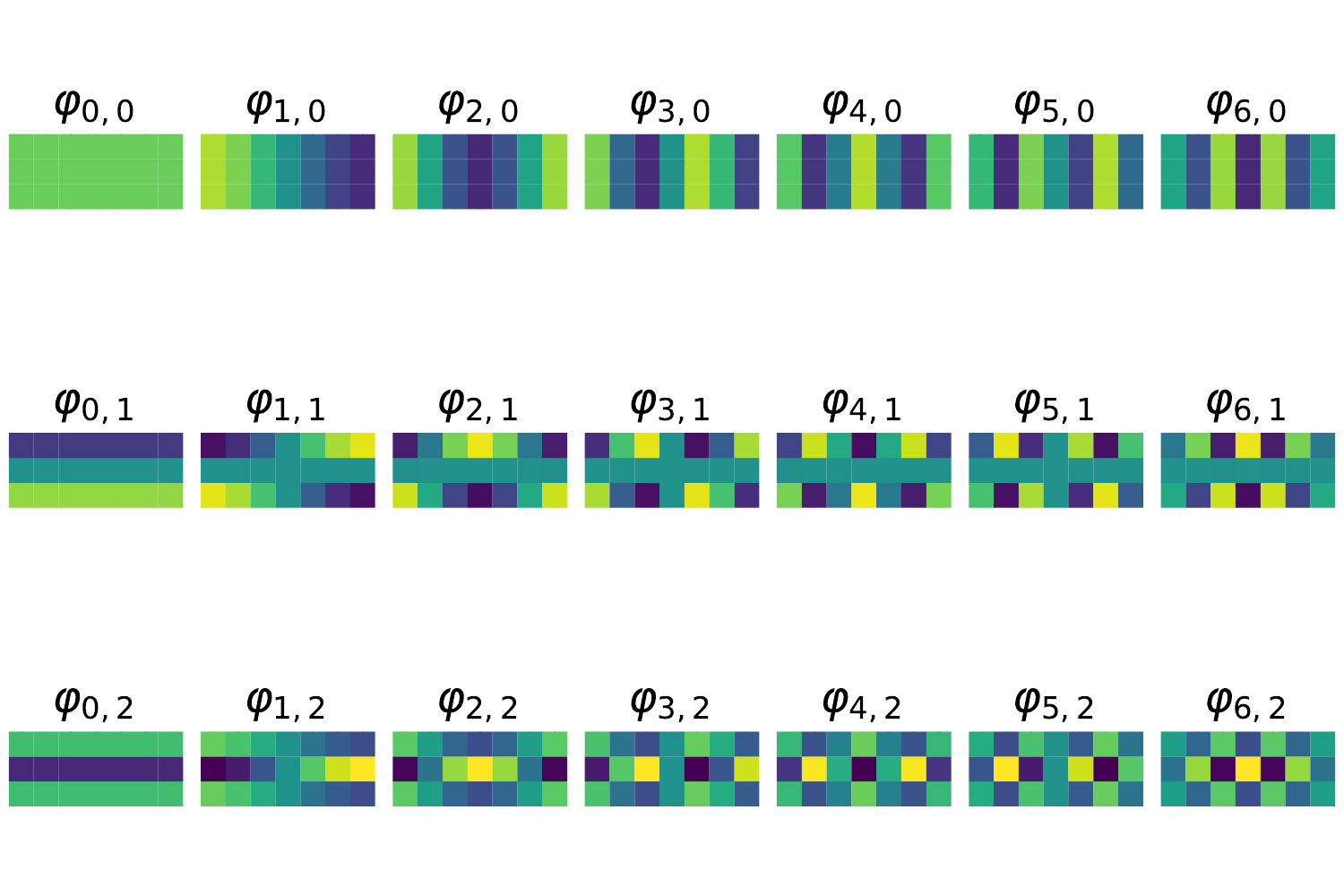}}
\vspace{-1.5em}
\caption{Laplacian eigenvectors of $P_7 \times P_3$ ordered in terms of
their natural horizontal/vertical frequencies.}
\label{fig:grid7x3dctorder}
\end{figure}

What we really want to do is to \emph{organize} those eigenvectors as shown in
Fig.~\ref{fig:grid7x3dctorder} instead of Fig.~\ref{fig:grid7x3evorder}.
\vspace{-1em}
Then a natural question is: how can we \emph{quantify the difference between the
  eigenvectors}? Note that the usual $\ell^2$-distance does not work since
$\| \bphi_i - \bphi_j \|_2 = \sqrt{2} \delta_{ij}$ where $\delta_{ij}$ is the
Kronecker delta.
Here, we propose to derive a natural distance between eigenvectors
using the ideas gained from the \emph{ramified optimal transport theory}
\cite{XIA-REVIEW}, i.e., we view
the cost to ``transport'' one eigenvector to another eigenvector
as the natural distance between two such eigenvectors.

\section{A Brief Review of Ramified Optimal Transport Theory}
\label{sec:ROT}
The Ramified Optimal Transport (ROT) theory \cite{XIA-REVIEW} is a branch of
more general optimal transport theory \cite{KOLOURI-ETAL-OMT-SPMAG}:
it studies transporting ``mass'' from one probability measure
$\mu^+$ to another $\mu^-$ along ramified transport paths with some specific
transport cost functional, and has been used to analyze various \emph{branching}
structures, e.g., trees; veins on a leaf; cardiovascular systems; river channel
networks; electrical grids; communication networks, to name a few.
For simplicity, we only consider the discrete case: two probability mass
functions (pmfs) in $\Rf^\ndim$.
Let $\f \define \sum_{i=1}^{k} f_i \delta_{\bx_i}$ and
$\bg \define \sum_{j=1}^{l} g_j \delta_{\by_j}$ be two such pmfs on the points
$\{\bx_i\}_{i=1}^k, \{\by_j\}_{j=1}^l \subset \Rf^\ndim$ with
$\sum_{i=1}^{k} f_i = \sum_{j=1}^{l} g_j = 1$.
Let $\Path(\f, \bg)$ be all possible transport paths from $\f$ to $\bg$
without cycles, i.e., each $G \in \Path(\f, \bg)$ is a weighted acyclic
directed graph with $\{ \bx_i \}_i \cup \{ \by_j \}_j \subset V(G)$,
whose edge~weights (> 0) satisfy \emph{the Kirchhoff law} at each node $v \in V(G)$:
\bdm \sum_{\substack{e \in E(G)\\ e^- = v}} w(e) = \sum_{\substack{e \in E(G)\\ e^+ = v}} w(e)
   + \begin{cases}
       f_i & \text{if $v=\bx_i$ for $\exists i \in 1:k$}\\
      -g_j & \text{if $v=\by_j$ for $\exists j \in 1:l$}\\
        0  & \text{otherwise}.
    \end{cases}
    \edm
Define the cost of a transport path $G \in \Path(\f, \bg)$ as:
\bdm
\bM_\alpha(G) \define \sum_{e \in E(G)} w(e)^\alpha \len(e),
\quad \alpha \in [0, 1].
\edm
Xia proved that the minimum transportation cost
$$ d_\alpha(\f,\bg) \define \min_{G \in \Path(\f,\bg)} \bM_\alpha(G) $$ is
a \emph{metric} on the space of pmfs and of homogeneous of degree $\alpha$; and
moreover he derived numerical algorithms to generate the $\alpha$-optimal path
for a given pair $(\f, \bg)$; see \cite{XIA-REVIEW} and the refereces therein
for further information.

\begin{hide}
Xia further derived the following facts among other things; see \cite{XIA-REVIEW} and the refereces therein for further information:
1) Number of branching nodes in $\Path(\f, \bg)$ can be bounded from above
by $k+l-2$;
2) The uniform lower bounds of minimum angle between any two edges in any
$\alpha$-optimal path in $\Path(\f,\bg)$;
3) The minimum transportation cost $\displaystyle d_\alpha(\f,\bg) \define
\min_{G \in \Path(\f,\bg)} \bM_\alpha(G)$ is a \emph{metric} on the space of
atomic measures of equal mass and of homogeneous of degree $\alpha$; and
4) Numerical algorithms to compute the $\alpha$-optimal path for a given
pair $(\f, \bg)$.
\end{hide}

\section{Our Proposed Method}
We propose the following:
\begin{alg}[LapEigPort]\ \\[-1em]
\label{alg:LapEigPort}
\begin{description}
\item[Step 0:] Convert $\bphi_i$ to a pmf $\bp_i$ over a graph $G$, $i=0:n-1$.
\item[Step 1:] Compute the cost to transport $\bp_i$ to $\bp_j$
  optimally, for all $i, j = 0:n-1$, which results in a ``distance'' matrix
  $D = (D_{ij}) \in \Rf_{\geq 0}^{n \times n}$.
\item[Step 2:] Use some embedding technique, e.g., \emph{Multidimensional Scaling} (MDS), to embed $D$ into a low dimensional Euclidean space $\Rf^{n_0}$, $n_0 \ll n$ (e.g., $n_0=2$ or $3$).
\item[Step 3:] Examine how eigenvectors are placed and organized in that embedding space
  $\Rf^{n_0}$.
\end{description}
\end{alg}
In Step~0, we set $\bp_i \define \(\phi_i^2[0], \ldots, \phi_i^2[n-1]\)^\transp$
in this paper; note that $\| \bphi_i \|_2 = 1$, $i=0:n-1$.
This is the most natural and obvious way to convert $\bphi_i$ to a pmf;
however, this may have some drawbacks and is certainly not the only way to
convert $\bphi_i$ to a pmf as we will discuss later.

The key step of Algorithm~\ref{alg:LapEigPort} is Step~1, which needs further
explanation.
Unlike the general ROT setting, a graph $G$ is fixed and given in our case.
On the other hand, our graph $G$ is \emph{undirected} while
the ROT requires \emph{directed} graphs.
Hence, we turn an undirected graph $G$ into the \emph{bidirected} graph
$\dtilde{G}$, a special type of directed multigraphs, i.e., each edge of $G$
joining two nodes, say, $i$ and $j$, becomes two directed edges
$(i,j)$ and $(j,i)$ in $\dtilde{G}$. 
To do so, we first compute the \emph{incidence matrix}
$Q = \left[ \bq_1 | \cdots | \bq_m \right] \in \Rf^{n \times m}$ of the
undirected graph $G=G(V,E)$ with $n=|V|$, $m=|E|$.
Here, $\bq_k$ represents the endpoints of the $k$th edge $e_k$: if $e_k$ joins nodes $i$ and
$j$, then $\bq_k[l]=1$ if $l=i$ or $j$; otherwise $\bq_k[l]=0$.
Then we orient each edge of $G$ in an artibrary manner to form
a directed graph $\tilde{G}$ and its incidence matrix
$\tilde{Q} = \left[ \widetilde{\bq}_1 | \cdots | \widetilde{\bq}_m \right] \in \Rf^{n \times m}$ where $\widetilde{\bq}_k$ is defined as follows:
if $e_k=(i,j)$ for some $i,j$, then $\widetilde{\bq}_k[i]=-1$,
$\widetilde{\bq}_k[j]=1$, and $\widetilde{\bq}_k[l]=0$ for $l \neq i, j$.
\begin{hide}
a directed graph $\tilde{G}$ whose incidence matrix $\tilde{Q}$ is, e.g.,
$\widetilde{\bq}_k[l] = -1$ if $l = i$; $=1$ if $l=j$; $=0$ otherwise.
 \bdm
 \widetilde{\bq}_k[l] = \begin{cases}
   -1 & \text{if $l = i$};\\
   1 & \text{if $l = j$}; \\
   0 & \text{otherwise.}
 \end{cases}
 \edm  
\end{hide}
Finally, form the bidirected graph $\dtilde{G}$ with $\dtilde{Q} \define \left[ \tilde{Q} \,\,\, | \, -\tilde{Q} \right] \in \Rf^{n \times 2m}$.

Given $\dtilde{Q}$, we solve the \emph{balance equation} that forces the
Kirchhoff law:
\begin{equation}
\label{eqn:balance}
  \dtilde{Q} \bw_{ij} = \bp_j - \bp_i, \quad
  \bw_{ij} \in \Rf^{2m}_{\geq 0}. 
\end{equation}
The weight vector $\bw_{ij}$ describes the transportation plan 
from $\bp_i$ to $\bp_j$. Note that $\bw_{ij}[k] = 0$, $\exists k \in [0, 2m)$
implies that the transportation plan does not use the corresponding
directed edge represented by the $k$th column of $\dtilde{Q}$.
Let $\dtilde{G}_{ij}$ be the bidirected graph $\dtilde{G}$
with these edge weights; then
$\dtilde{G}_{ij} \in \Path(\bp_i, \bp_j)$.
We caution here that Eq.~\eqref{eqn:balance} may have multiple solutions.
Hence, we propose to solve the following \emph{Linear Programming} (LP) problem (see, e.g., \cite{PAPA-STEIGLITZ}):
\begin{equation}
\label{eqn:l1min}
\min_{\bw_{ij} \in \Rf^{2m}} \| \bw_{ij} \|_1
\enskip \text{subj.\ to:}
\begin{cases}
\dtilde{Q} \bw_{ij} = \bp_j - \bp_i ; \\
\bw_{ij}[l] \geq 0, l=0:2m-1,
\end{cases}
\end{equation}
to obtain one of the \emph{sparse} solutions of Eq.~\eqref{eqn:balance},
which turned out to be much faster to compute than using the
nonnegative least squares (NNLS) solver \cite[Chap.~23]{LAWSON-HANSON},
which tends to generate denser solutions. 

Finally we fill the distance matrix $D=(D_{ij})$ by
\bdm
D_{ij} = \bM_\alpha( \dtilde{G}_{ij} ) = \sum_{e \in E(\dtilde{G}_{ij})} w_{ij}(e)^\alpha \len(e), \quad \alpha \in [0, 1].
\edm
Note that Eq.~\eqref{eqn:l1min} may still have multiple solutions,
and currently we are \emph{not} examining all possible solutions
of Eq.~\eqref{eqn:l1min} to search
$\arg \min_{\dtilde{G}_{ij} \in \Path(\bp_i, \bp_j)} \bM_\alpha(\dtilde{G}_{ij})$.
We plan to investigate how to handle such multiple solutions in LP, e.g.,
the method proposed in \cite{DYER-PROLL}.

\section{Numerical Results}
In this section, we will demonstrate the effectiveness of our proposed method
using the same 2D lattice graph we discussed earlier as well as a dendritic tree
of a retinal ganglion cell (RGC) of a mouse.
We note that we used the \texttt{JuMP} optimization package \cite{JUMP} written
in \texttt{Julia} \cite{JULIA} in order to solve the optimization problem
Eq.~\eqref{eqn:l1min}.

\vspace{-1em}
\subsection{The 2D lattice graph $P_7 \times P_3$}
\begin{figure}[h]
\centering\includegraphics[width=0.48\textwidth]{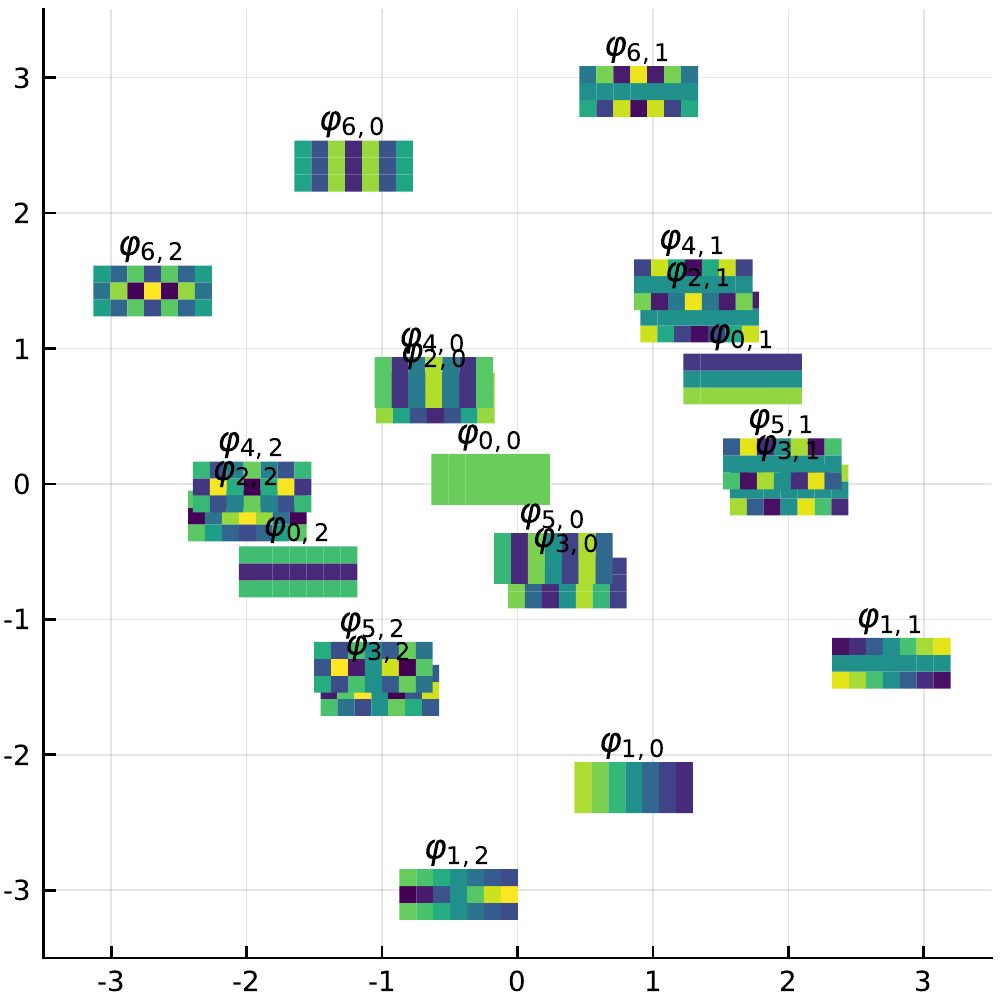}
\caption{Embedding of the Laplacian eigenvectors of $P_7 \times P_3$ into
  $\Rf^2$ using Algorithm~\ref{alg:LapEigPort} with $\alpha=0.5$.}
\label{fig:grid7x3map}
\end{figure}

Fig.~\ref{fig:grid7x3map} shows the embedding of these 21 eigenvectors into
$\Rf^2$ using Algorithm~\ref{alg:LapEigPort} where we set $\alpha=0.5$ and
used the classical MDS \cite[Chap.~12]{BORG-GROENEN} for embedding.
Of course, in general, when a graph is given, we cannot assume the best
embedding dimension $n_0$ a priori. Here we simply used $n_0=2$ because the
top two eigenvalues of the Gram matrix of the configurations (i.e., the outputs
of the MDS) were more than twice the third eigenvalue.
Fig.~\ref{fig:grid7x3map} clearly reveals the two-dimensional
organization of the eigenvectors, and is somewhat similar to a rotated version
of Fig.~\ref{fig:grid7x3dctorder}, but the details are different.
For example, the eigenvectors with even and odd oscillations are 
mapped in a symmetric manner around the ``DC'' vector $\bvphi_{0,0}$ located
at the center of Fig.~\ref{fig:grid7x3map}. This symmetry is due to the use of
the energy of the eigenvector components in Step~0 of
Algorithm~\ref{alg:LapEigPort}: from Eq.~\eqref{eqn:eigen1},
we have $\phi_{k;n}^2[x] + \phi_{n-k;n}^2[x] \equiv a_{k;n}^2$,
$k=1:n-1$, $x=0:n-1$.

We also examined the case of $\alpha=1$, which gave us more congested embedding
around $\bvphi_{0,0}$.

\vspace{-1em}
\subsection{The dendritic tree of an RGC of a mouse}
\begin{figure}[htb]
\begin{minipage}[b]{0.5\linewidth}
  \centering
  \centerline{
  \includegraphics[width=\textwidth]{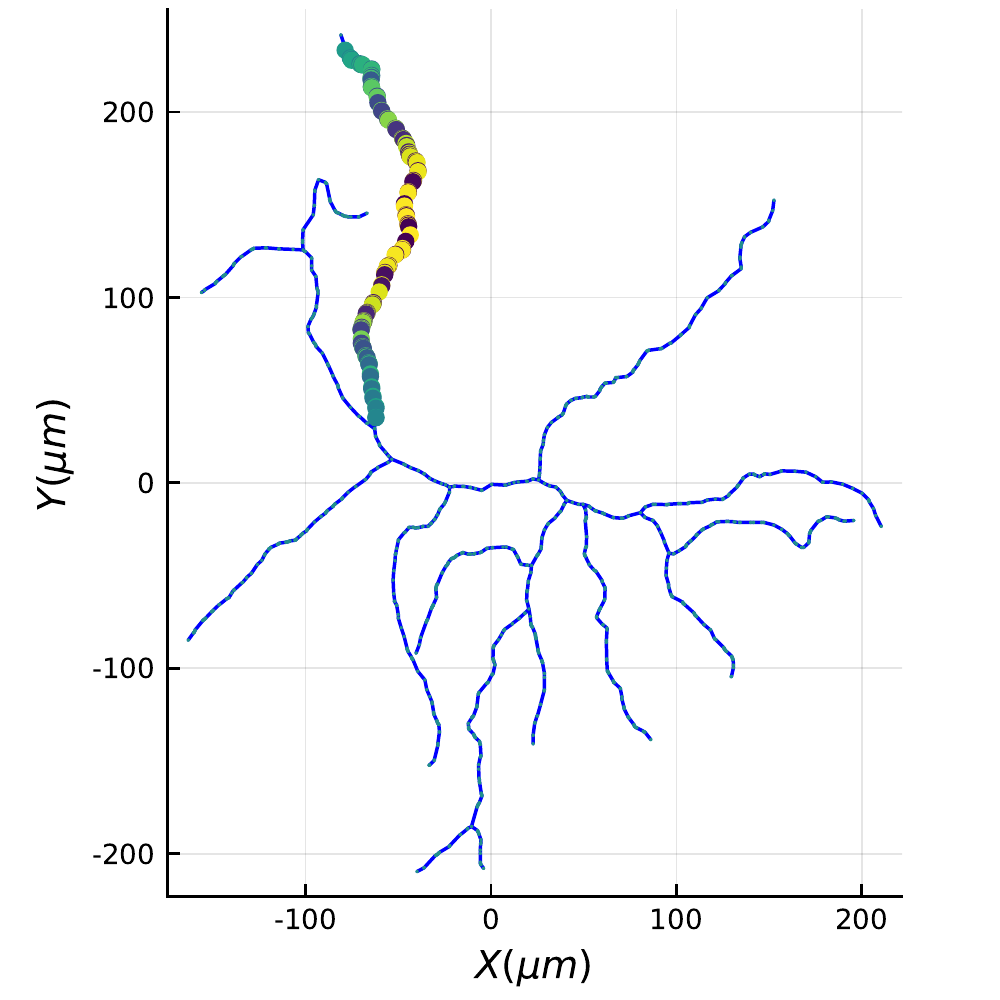}
  }
  \centerline{(a) $\bphi_{1141}$ with $\lambda_{1141}=3.9994$}\medskip
\end{minipage}
\hfill
\begin{minipage}[b]{0.5\linewidth}
  \centering
  \centerline{
  \includegraphics[width=\textwidth]{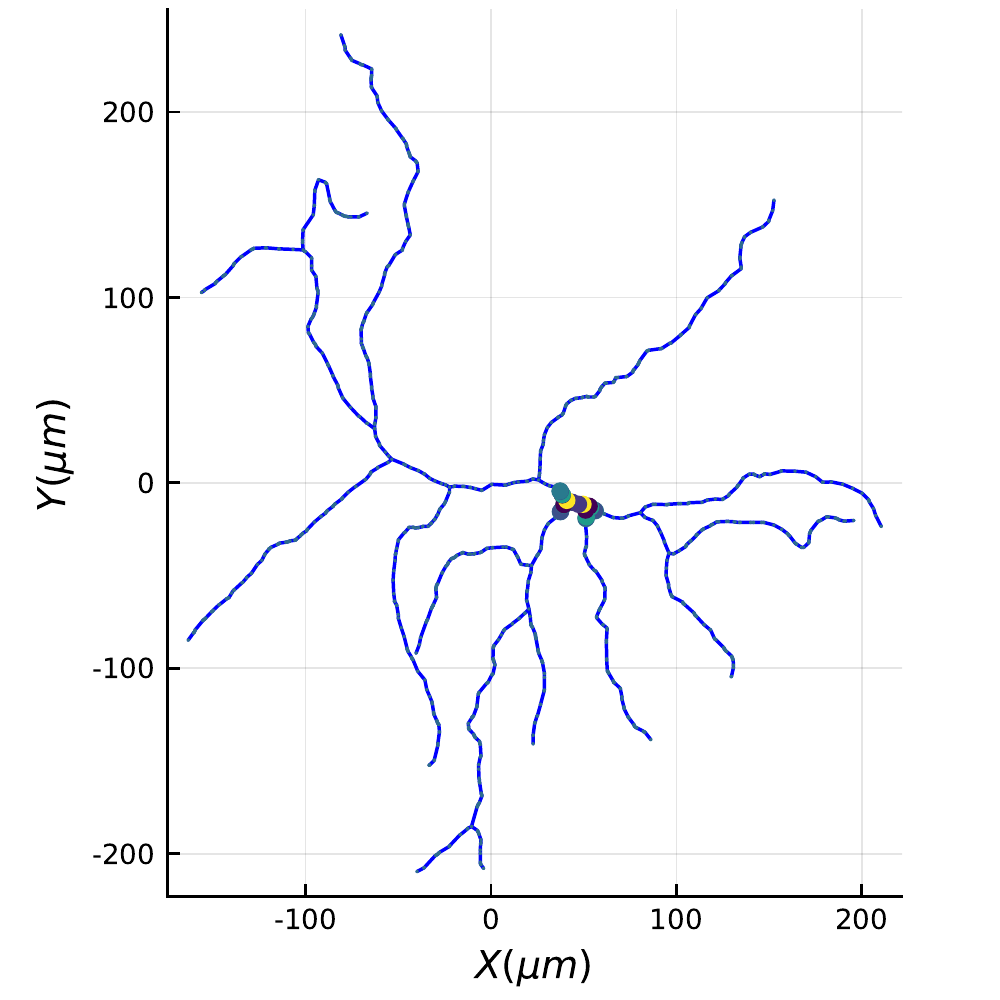}
}
  \centerline{(b) $\bphi_{1142}$ with $\lambda_{1142}=4.3829$}\medskip
\end{minipage}
\vspace{-3em}
\caption{Laplacian eigenvectors before and after the phase transition.
  The yellow and purple circles indicate larger positive and negative 
  components while those on thin blue lines are negligible.
  This tree has $n=1154$ nodes in total.}
\label{fig:RGCphasetrans}
\end{figure}

\begin{figure}[h]
\centering\includegraphics[width=0.5\textwidth]{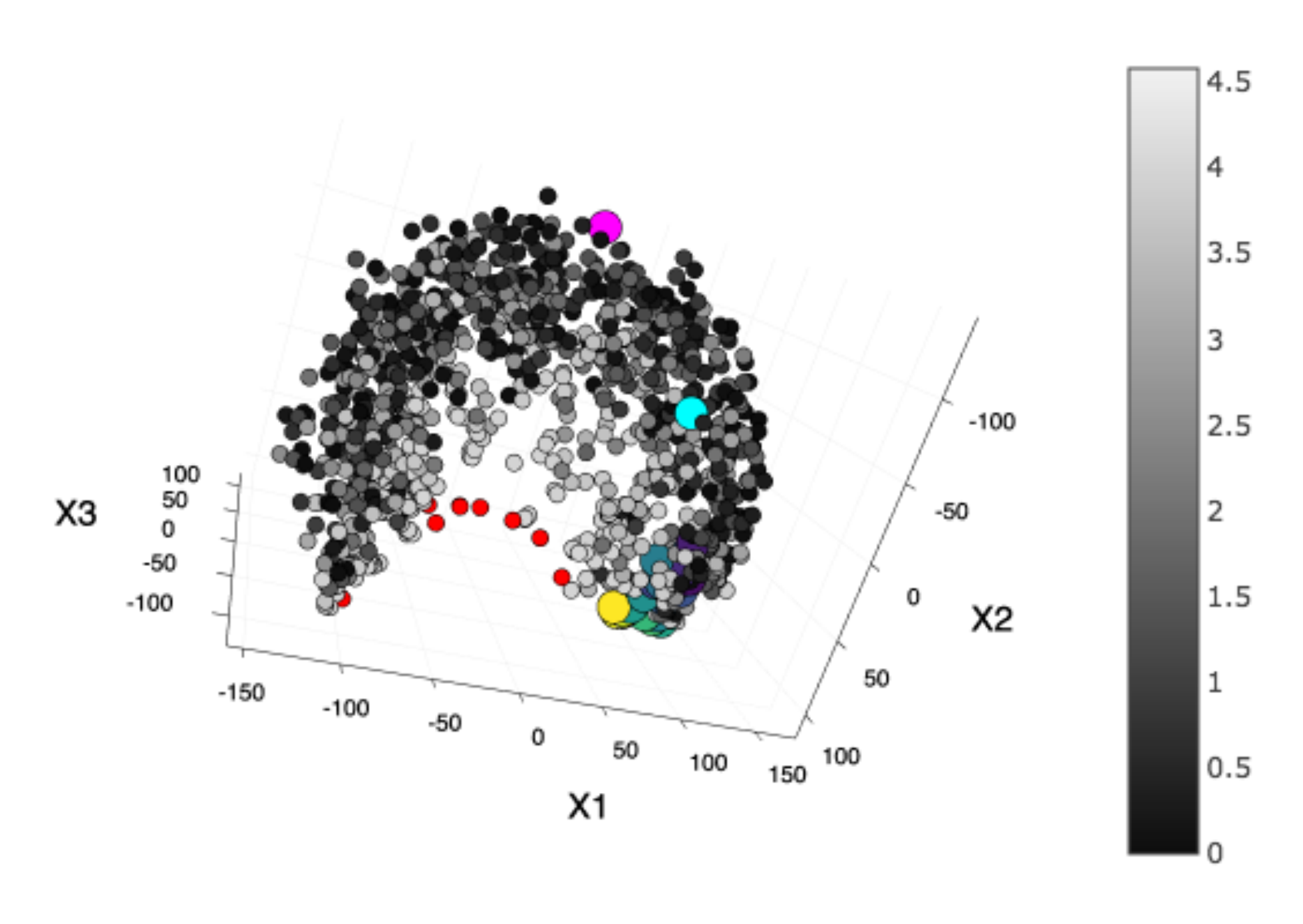}
\vspace{-3em}
\caption{Embedding of the Laplacian eigenvectors of the RGC tree into
  $\Rf^3$ using Algorithm~\ref{alg:LapEigPort} with $\alpha=0.5$.}
\label{fig:RGCR3}
\end{figure}

For the details on the data acquisition and the conversion of dendritic trees of
RGCs to graphs (in fact, literally ``trees''), see \cite{SAITO-WOEI-DENDRITES}.
As discussed in \cite{SAITO-WOEI-KOKYUROKU, WHY4-LAA} in greater detail,
the graph Laplacian eigenvectors of these trees exhibit a very peculiar
\emph{phase transition phenomenon}: $\bphi_k$'s with $\lambda_k < 4$ oscillate
\emph{semi-globally} while those with $\lambda_k > 4$ are \emph{concentrated}
around junctions (i.e., nodes whose degrees are greater than 2).
Fig.~\ref{fig:RGCphasetrans} shows 2D projections of two of the eigenvectors
of a particular dendritic tree demonstrating this phenomenon.
Fig.~\ref{fig:RGCR3} shows the embedding of the Laplacian eigenvectors of this
dendritic tree into $\Rf^3$ using Algorithm~\ref{alg:LapEigPort}
with $\alpha=0.5$.
We used $n_0=3$ because of the size difference between the top three
eigenvalues and the rest. 
These 3D point clouds form a upside-down bowl-like shape with three legs
(or a shape partly joining two croissants).
The magenta circle indicates the DC vector, $\bphi_0$, while the cyan circle
is the Fiedler vector, $\bphi_1$.
The red circles located around the bottom are the eigenvectors $\bphi_{k}$, $k > 1141$, i.e., those localized around the junctions. The larger colored circles
located in the lower right region are the eigenvectors supported on
one of the branch indicated in Fig.~\ref{fig:RGCphasetrans}(a).
The gray scale colors of the other points indicate the magnitude of the
corresponding eigenvalues.
\vspace{-2.25em}

\vspace{-1.5em}
\section{Discussion}
Although Algorithm~\ref{alg:LapEigPort} allows us to see interpretable
organizations of the graph Laplacian eigenvectors,
there remain several questions that need to be
answered: 1) How can we choose optimal value of $\alpha \in [0,1]$?
2) Why the Fiedler vectors tend to be mapped far from the DC vectors?
3) How can we find the true cost-minimizing transportation path among
$\Path(\bp_i, \bp_j)$?
4) What other ways to turn $\bphi_i$ into $\bp_i$ should we consider?
\begin{hide}
Generating $\bphi_i^2$ is not the only way to turn $\bphi_i$ into a pmf.
Other possibilities include: 1) an $\ell^1$-normalized pmf: $\bphi_i^1 \define \( \left| \phi_i[0] \right|, \dots, \left| \phi_i[n-1] \right| \)^\transp / \| \bphi_i \|_1$ ; 2) an inflated pmf (via adding a constant followed by normalization):
\bdm
\widetilde{\bphi}_i \define \begin{cases}
  \bphi_0^1 & \text{if $i=0$;} \\
  \frac{\bphi_i - c_\mathrm{min} \cdot \bone_n}{\left\| \bphi_i - c_\mathrm{min} \cdot \bone_n \right\|_1} & \text{if $i \neq 0$},
\end{cases}
\edm
where $\displaystyle c_\mathrm{min} \define \min_{0 < i < n; \, 0 \leq l < n} \bphi_i[l] < 0$.
\end{hide}
5) How can we develop the true Littlewood-Paley theory on graphs and most
``natural'' wavelets on graphs once the above embedding is done?

\vfill
\pagebreak


\end{document}

%% file: amsmath_def.tex
\chardef\bslash=`\\ 

\makeatletter
\def\verbatim{\interlinepenalty\@M \@verbatim
  \leftskip\@totalleftmargin\advance\leftskip2pc
  \frenchspacing\@vobeyspaces \@xverbatim}
\makeatother
\input{amsmath_sym}

\theoremstyle{plain}   
\newtheorem{thm}{Theorem}[section]   
\newtheorem{alg}[thm]{Algorithm}     

\theoremstyle{definition}



%% file: amsmath_sym.tex
\newcommand{\bdm}{\begin{displaymath}}
\newcommand{\edm}{\end{displaymath}}

\newcommand{\half}{\frac{1}{2}}

\renewcommand{\(}{\left(}
\renewcommand{\)}{\right)}
\newcommand{\N}{{\mathbb N}}

\newcommand{\Rf}{{\mathbb R}}

\newcommand{\define}{\, := \,}

\newcommand{\udots}{
  \mathinner {\mkern 1mu\raise 1pt \vbox {\kern 7pt \hbox {.}}\mkern 2mu
  \raise 4pt \hbox {.}\mkern 2mu\raise 7pt \hbox {.}\mkern 1mu}}

\newcommand{\f}{{\boldsymbol f}} 
\newcommand{\bg}{{\boldsymbol g}}

\newcommand{\bp}{{\boldsymbol p}}
\newcommand{\bq}{{\boldsymbol q}}

\newcommand{\bw}{{\boldsymbol w}}
\newcommand{\bx}{{\boldsymbol x}}
\newcommand{\by}{{\boldsymbol y}}

\newcommand{\bM}{{\boldsymbol M}}

\newcommand{\bphi}{{\boldsymbol \phi}}
\newcommand{\bvphi}{{\boldsymbol \varphi}}

%% file: symbols.tex
\newcommand{\ndim}{n}

\newcommand{\transp}{{\scriptscriptstyle{\mathsf{T}}}}

\newcommand{\bone}{{\boldsymbol 1}}



%% file: abs2.tex
When attempting to develop wavelet transforms for graphs and networks,
some researchers have used graph Laplacian eigenvalues and eigenvectors in place
of the frequencies and complex exponentials in the Fourier theory for regular
lattices in the Euclidean domains. This viewpoint, however, has a fundamental
flaw: on a general graph, the Laplacian eigenvalues \emph{cannot} be interpreted
as the \emph{frequencies} of the corresponding eigenvectors. 
In this paper, we discuss this important problem further and propose a new
method to organize those eigenvectors by defining and measuring ``natural''
distances between eigenvectors using the Ramified Optimal Transport Theory 
followed by embedding them into a low-dimensional Euclidean domain.
We demonstrate its effectiveness using a synthetic graph as well as a dendritic
tree of a retinal ganglion cell of a mouse.